\documentclass[11pt]{article}
\usepackage{amsmath,amssymb,amsfonts}
\usepackage[T2A]{fontenc}

\newtheorem{theorem}{Theorem}

\begin{document}

\title{Surfaces via spinors and soliton equations
\thanks{This work was partially supported by the Mathematical Center in Akademgorodok under
the agreement No. 075-15-2019-1675 with the Ministry of Science and Higher Education
of the Russian Federation.}}
\author{Iskander A. Taimanov
\thanks{Sobolev Institute of Mathematics, 630090 Novosibirsk, Russia, and Department of Mathematics
and Mechanics, Novosibirsk State University, 630090 Novosibirsk, Russia, e-mail:
taimanov@math.nsc.ru.}}
\date{}
\maketitle

\abstract{
This article surveys the Weierstrass representation of surfaces in the three- and four-dimensional
spaces, with an emphasis on its relation to the Willmore functional. We also
describe an application of this representation to constructing a new type of solutions to the
Davey--Stewartson II equation. They have regular initial data, gain one-point singularities
at certain moments of time, and extend to smooth solutions for the remaining times.}

\section{The Weierstrass (spinor) representation of surfaces in the three-space}

The Weierstrass representation for minimal surfaces in the three-space
is as follows: for any pair of holomorphic functions
$\psi_1$ and $\bar{\psi}_2$ defined in a domain $U \subset {\mathbb C}$ in a complex plane
the formulae
$$
x^1(P) = \frac{i}{2} \int [(\psi_1^2 + \bar{\psi}_2^2)dz + (\bar{\psi}_1^2 + \psi_2^2)d\bar{z}] + x^1(P_0),
$$
\begin{equation}
\label{weier}
x^2(P) = \frac{1}{2} \int [(-\psi_1^2 + \bar{\psi}_2^2)dz + (-\bar{\psi}_1^2 + \psi_2^2)d\bar{z}] + x^2(P_0),
\end{equation}
$$
x^3(P) = \int[\psi_1 \bar{\psi}_2 dz + \bar{\psi}_1\psi_2 d\bar{z}] + x^3(P_0)
$$
determine a minimal surface in ${\mathbb R}^3$. Here we assume that $U$ is simply-connected or the integrals over cycles in $U$ vanish and the integrals are taken along a path from a fixed point $P_0 \in U$ to $P$. Moreover  every minimal surface admits such a representation. Weierstrass used another data namely
$f = \bar{\psi}_2^2$ and $g = \frac{\psi_1}{\bar{\psi}_2}$.  However for the generalization of this representation it is worth to consider $\psi_1$ and $\psi_2$ and treat this pair as a solution of the Dirac equation
\begin{equation}
\label{diraceq}
D \psi = 0, \ \ \ \psi = \left(\begin{array}{c} \psi_1 \\ \psi_2 \end{array}\right),
\end{equation}
for a two-dimensional Dirac operator of the form
\begin{equation}
\label{diracop}
D = \left(\begin{array}{cc} 0 & \partial \\ -\bar{\partial} & 0 \end{array}\right) + 
\left(\begin{array}{cc} U & 0 \\ 0 & U \end{array}\right), \ \ U = \bar{U},
\end{equation}
where a real-valued potential $U$ vanishes for minimal surfaces.
Now the Weierstrass representation generalizes as follows

\begin{theorem}[\cite{K1}]
For every solution $\psi$ of (\ref{diraceq}) the formulae (\ref{weier}) define a surface in ${\mathbb R}^3$
for which $z$ is a conformal parameter, induced
metric takes the form
$$
ds^2 = e^{2\alpha}dz d\bar{z}, \ \ \ \ e^\alpha = |\psi_1|^2 +
|\psi_2|^2
$$
and the potential $U$ of the Dirac operator equals to
$$
U = \frac{He^\alpha}{2},
$$
where $H$ is the mean curvature.
\end{theorem} 

\begin{theorem}[\cite{T1}]
Every surface in ${\mathbb R}^3$ (with a fixed conformal parameter $z$ on it) 
admits such a representation even globally. Therewith $\psi$ is a section of
a spinor bundle over the surface, the form $U^2 dx \wedge dy$ is globally defined and its integral over the surface is proportional to the Willmore functional
$$
{\cal W} = \int H^2 d\mu = 4 \int U^2 dx \wedge dy,
$$
where $d\mu$ is the induced area form of the surface.
\end{theorem}

Hence being considered for the Dirac operators with general real-valued potentials the formulae (\ref{weier}) define the {\it Weierstrass (spinor) representation} of general surfaces in ${\mathbb R}^3$.

Theorem 1 was derived from the similar formulae in the book by Eisenhart \cite[Problem 35.4]{Eisenhart} where instead of (\ref{diraceq}) the following condition is used:
$$
L \psi_1 = L \bar{\psi}_2 = 0, \ \ \ L = \partial \bar{\partial} - \frac{\partial \log U}{U}\bar{\partial} + U^2.
$$
$D$ naturally arises as the ``square root'' of the Schr\"odinger operator $L$. The representation basing on the Dirac operator has much more opportunities because its potential has no singularities and the operator has good spectral properties. In the advanced problems of his textbook Eisenhart frequently proposed to prove results from various articles and we cannot exclude that these formulae are traced to some earlier publication. It appears that this local representation is equivalent to another one derived in \cite{Kenmotsu} where  the Dirac operator was not used either.  

In \cite{K1} the Weierstrass representation was used for 
introducing the deformations of surfaces admitting such a representation. 
The operator $D$ generates a hierarchy of solution equations of the form
$$
\frac{\partial D}{\partial t_n} = [D, A_n] - B_n D,
$$
where $A_n$ and $B_n$ are matrix differential operators such that
the principal term of $A_n$ takes the form
$$
A_n = \left(\begin{array}{cc} \partial^{2n+1} +
\bar{\partial}^{2n+1} & 0 \\ 0 & \partial^{2n+1} +
\bar{\partial}^{2n+1}
\end{array}\right) + \dots\ .
$$
This evolution preserves 
the zero energy level of $D$ deforming the corresponding
eigenfunctions
\begin {equation}
\label{psievol}
\frac{\partial \psi}{\partial t} + A \psi = 0
 \end{equation} 
and $D\psi_0 = 0$ for the initial data $\psi_0=\psi|_{t=t_0}$,
then $D \psi = 0$ for all $t \geq t_0$.

For $n=1$ we have the modified Novikov--Veselov (mNV) equation 
\cite{Bogdanov}
\begin{equation} 
\label{mnv} 
U_t = \left(U_{zzz} + 3 U_z V + \frac{3}{2} U V_z \right) +
\left(U_{\bar{z}\bar{z}\bar{z}} + 3 U_{\bar{z}} \bar{V} +
\frac{3}{2} U \bar{V}_{\bar{z}}\right) 
\end{equation}
where
\begin{equation}
\label{mnv-cons} V_{\bar{z}} = (U^2)_z.
\end{equation}
In the case when $U\vert_{t=0}$
depends only on $x$ we have $U = U(x,t)$ and the mNV equation
reduces to the modified Korteweg--de Vries equation 
$U_t = \frac{1}{4} U_{xxx} + 6U_x U^2$ (here $V =
U^2$). In the same manner the original Novikov--Veselov equation 
$$
U_t = U_{zzz}+ U_{\bar{z}\bar{z}\bar{z}} + (V U)_z + (\bar{V}
U)_{\bar{z}}, \qquad V_{\bar{z}} = 3  U_z
$$
generalizes the Korteweg--de Vries equation.

The mNV deformation introduced in \cite{K1} is as follows:
let a surface be induced by $\psi$ via (\ref{weier}) and consider solutions 
$U$ and $\psi$ of the mNV equation and (\ref{psievol}) with given initial data. 
Then for any moment of time we have a spinor $\psi$ that determines the deformed surface.
In fact, we have infinitely many deformations defined up to translations by $(x^1(P_0,t), x^2(P_0,t),x^3(P_0,t))$. 
This is some family of the {\it mNV deformations} of the surface. 

\begin{theorem}[\cite{T1}]
The mNV deformations evolve tori into tori and preserve their conformal classes and the values of
the Willmore functional.
\end{theorem}

Theorems 2 and 3 hint on the relation of this representation to the Willmore functional and in Section 2 
we briefly expose how it was applied to studying the conformal geometry of surfaces. In Section 4, in difference with Section 2 where analysis was applied to geometry, we discuss
the recent applications of geometry to analysis. We show how to construct exact solutions
to the Davey--Stewartson II equation. Therewith, geometry of surfaces helps to find a new scenario for creating 
singularities of solutions with regular initial data.

It would be interesting to apply the Weierstrass representation to other problems of surface theory
(bending, existence of umbilics, etc.). In particular, if some conjecture appears false, then methods of integrable systems can help in constructing an explicit counterexample (see, for instance, \cite{Abresch}).

\section{Spectral characteristics of $D$ and conformal geo\-metry of surfaces}

The Willmore conjecture which reads that the minimum of the Willmore functional among tori in 
${\mathbb R}^3$ is attained at the Clifford torus was proved in \cite{MN} by means of geometric measure theory and calculus of variations.

In the mid-1990s we proposed an approach to proving it that was based on Theorem 3 and integrable systems theory.
This approach was not implemented but we think it is worth to be briefly exposed here. 

It was conjectured in \cite{T1} that 

{\sl a torus nonstationary (with respect to the mNV flow and up to translations) cannot be a local minimum of the Willmore functional}. 

\noindent
Otherwise, by Theorem 3, the minimum of the Willmore functional has to contain an infinite family of tori invariant under the mNV flow and that is very unlikely. By the general philosophy of integrable systems, the stationary solution to the mNV equation has the simplest possible spectral curve \cite{T2}. 

For two-dimensional differential operators with periodic coefficients the spectral curve 
(on the zero level energy) parameterizes its Floquet eigenfunctions \cite{DKN}.
In our case a Floquet eigenfunction $\psi$ of the operator $D$ with the
eigenvalue (or the energy) $E$ is a formal solution to the equation
$$
D \psi = E \psi
$$
which satisfies the periodicity conditions
$$
\psi(q+\gamma_j) = e^{2\pi i (k,\gamma_j)} \psi(z,\bar{z}), \ \ j=1,2,
$$
where $\gamma_1$ and $\gamma_2$ generate the lattice of periods $\Lambda$ of the potential $U$ and
$(k,\gamma) = k_1 \gamma_j^1 + k_2 \gamma_j^2$ is the inner product. 
The quantities $k_1,k_2 \in {\mathbb C}$ are called the quasimomenta of $\psi$ and 
$\mu(\gamma_j) =e^{2\pi i (k,\gamma_j)}$ are Floquet multipliers.
All possible triples $(k_1,k_2,E)$ for which Floquet functions exist form an analytic subset $Q(U)$ 
in ${\mathbb C}^3$ invariant under the dual lattice 
$\Lambda^\ast \subset {\mathbb R}^2 \subset {\mathbb C}^2$ acting on the quasimomenta. 
We proved that for the two-dimensional operators $\Delta +U$ and $\partial_y - \partial^2_x + U$ in 1985. However this paper was unpublished, although referred in \cite{K89} and was exposed in \cite{T06}.
Now we define the spectral curve as the complex curve
$$
\Gamma = (Q \cap \{E=0\})/\Lambda^\ast
$$
and consider it up to biholomorphic equivalence making the definition independent
on the choice of a basis for $\Lambda$. The curve is an invariant of the mNV flow, it is naturally 
completed by a couple of points at infinity which compactify it in the case of finite genus. The Floquet functions are glued into a meromorphic section over $\Gamma$. The above rough definition must be  detailed for singular spectral curves. In general the space of Floquet functions corresponding to a point from $\Gamma$ is one-dimensional and the multiple points have to be normalized in such a manner that for the resulted curve $\Gamma_\psi$  to every point there corresponds a one-dimensional space, there is a meromorphic section $\psi$ of this bundle and every Floquet function is a linear combination of sections at different points (see the definition of $\Gamma_\psi$ in \cite{T06}). The spinor $\psi$ generating a torus via (\ref{weier}) has the Floquet multipliers equal to $\pm 1$.

The spectral curve defined for $D$ is a particular case of the general spectral curves
which play fundamental role in integrable systems. They are the first integrals of the system 
(that was first showed for the Korteweg-de Vries equation in \cite{Novikov}). The particular case of them
are the spectral curves of constant mean curvature tori which are always of finite genus \cite{Hitchin,PS}.
In general this spectral curve is of infinite genus. For finite genera cases solutions to the integrable systems are expressed in terms of theta functions on spectral curves. In our case all Floquet functions are reconstructed from certain data related to $\Gamma_\psi$ and the value of the Willmore functional is also determined by them \cite{T2}. We conjectured that 

{\sl for tori in ${\mathbb R}^3$ the curve $\Gamma$, i.e., the set of the multipliers $\mu(\gamma_j)$,  is conformally invariant (as well as the Willmore functional).} 

\noindent
Since this is evident for translations and rotations, that was left to prove for the Mobius inversion and that was confirmed in \cite{GS}.

For the Clifford torus parameterized by $x,y$ such that $0 \leq x,y\leq 2\pi$ the potential $U$
of its Weierstrass representation is 
$$
U(x) = \frac{\sin x}{2\sqrt{2}(\sin x - \sqrt{2})}
$$
and its spectral curve $\Gamma_\psi$ is ${\mathbb C}P^1$ with two pairs of glued points.

For differential operators on surfaces of higher genera the analog of Floquet--Bloch theory is unknown. It would be interesting to find it, if it exists, at least for the Dirac operator $D$.

For spheres there are no analogs of the Floquet functions and the zero energy level of $D$ just consists of  the kernel $\mathrm{Ker}D$. We notice that there is an antiinvolution
\begin{equation}
\label{symmetry}
\left(\begin{array}{c} \psi_1 \\ \psi_2 \end{array}\right) \stackrel{\sigma}{\longrightarrow}
\left(\begin{array}{c} - \bar{\psi}_2 \\ \bar{\psi}_1 \end{array}\right), \ \ \sigma^2 = -1,
\end{equation}
acting on $\mathrm{Ker}D$. This implies that the dimension of the kernel over ${\mathbb C}$ is 
always even. 

We say that a sphere in ${\mathbb R}^3$ admits a spinor representation with a one-dimensional potential
if after removing a certain pair of points we obtain the cylinder ${\mathbb R} \times S^1$ for which the potential of the representation depends on $x$ only; i.e.,  $U=U(x)$. These are, for instance, spheres of revolution. 
By using the inverse scattering transform of one-dimensional Dirac operators on the line we proved 

\begin{theorem}[\cite{T21}]
For spheres with a one-dimensional potential we have
\begin{equation}
\label{conj1}
{\cal W} = 4\int U^2 dx \wedge dy \geq  4\pi N^2,
\end{equation}
where $\dim_{\mathbb C} \mathrm{Ker}D = 2N$
and the equalities are achieved at the soliton potentials
$$
U_N(x) = \frac{N}{2\cosh x}.
$$
\end{theorem}

We call the spheres that correspond to these potentials soliton spheres and it appears that they have very interesting geometrical properties \cite{BP}. In \cite{T2} we conjectured that

{\sl inequality (\ref{conj1}) holds for all spheres}. 

Soon after the preprint of \cite{T2} appeared Friedrich showed that this conjecture implies the following statement:

{\sl Given an eigenvalue $\lambda$ of the Dirac operator $D$  on a two-dimensional spin-manifold homeomorphic to the two-sphere, 
\begin{equation}
\label{friedrich}
\lambda^2 \mathrm{Area}(M) \geq \pi m^2(\lambda) 
\end{equation}
where $m(\lambda)$ is the multiplicity of $\lambda$. }

For $m(\lambda)=2$ inequality (\ref{friedrich}) was already proved by B\"ar \cite{Bar}.

The arguments by Friedrich were as follows. On a spin-manifold of dimension $2$ with the metric $e^{2\alpha}  dz d\bar{z}$ the Dirac operator (on the spin-manifold) takes the form
$$
D = 2 e^{-3\alpha/2} \left(\begin{array}{cc} 0 & \partial \\ -\bar{\partial} & 0 \end{array}\right)e^{\alpha/2},
$$
and the equation
$$
D\varphi = \lambda \varphi
$$
is rewritten as
$$
\left[\left(\begin{array}{cc} 0 & \partial \\ -\bar{\partial} & 0 \end{array}\right) - \frac{\lambda e^\alpha}{2}\right]\psi = 0
$$
where $\psi = e^{\alpha/2}\varphi$, and if $\lambda$ is constant, then (\ref{conj1}) implies (\ref{friedrich}).
Moreover, if $\lambda = H$, then  this is exactly the Dirac equation
(\ref{diraceq}) (the sign of the mean curvature can be changed without any loss) and, since $e^\alpha = |\psi|^2$, we have $|\varphi|=1$. Therefore the Weierstrass representation is rewritten in terms of solutions 
of the Dirac equation
$$
D \varphi = H \varphi
$$
of constant length: $|\varphi|=1$ \cite[Theorem 13]{Friedrich}.

This embedding of the Weierstrass representation into the general framework of Dirac operators on spin-manifolds appears very fruitful: it led to its generalization, the spinorial representation of immersions of manifolds, which are not necessarily  two-dimensional, into certain homogeneous spaces (see \cite{Bayard} and references therein).

The Weierstrass representation for surfaces in ${\cal R}^3$ was generalized for surfaces in three-dimensional Lie groups with left-invariant metrics  in \cite{BT}. It helped to establish some facts on constant mean curvature surfaces in these groups.

It would be interesting, at least as a test problem, to find a discretization of the Weierstrass representation
by means of discrete complex analysis. In \cite{Zakharov} that was done for the generalizations of the representation for 
time-like surfaces in ${\mathbb R}^{2,1}, {\mathbb R}^{3,1}$ and ${\mathbb R}^{2,2}$. But in these cases complex analysis is not involved because the principal term of the Dirac operator $D$ has the form
$\left(\begin{array}{cc} 0 & \partial_\xi \\ \partial_\eta & 0 \end{array}\right)$ where $\xi$ and $\eta$ are isotropic coordinates.

The conjectured inequality (\ref{conj1}) was finally proved in \cite{FLPP} with its generalizations for surfaces of higher genera:

\begin{theorem}[\cite{FLPP}]
\label{pluck}
For a closed oriented surface of genus $g$ immersed into ${\mathbb R}^3$
via (\ref{weier}) and (\ref{diraceq}), we have
$$
\int U^2 dx \wedge dy \geq
\begin{cases}
\pi N^2 & \text{for $g=0$} \\
\begin{cases}
\frac{\pi N^2}{4} & \text{for $N$ even} \\
\frac{\pi (N^2-1)}{4} & \text{for $N$ odd}
\end{cases} 
& \text{for $g=1$} \\
\frac{\pi}{4g} \left( N^2 - g^2 \right) & \text{for $g > 1$},
\end{cases}
$$
where $\dim_{\mathbb C} \mathrm{Ker}D = 2N$.
\end{theorem}

\section{Surfaces in the four-space and the Davey--Ste\-wart\-son equation}

Theorem 7 was derived from the Pl\"ucker formula in the quaternionic algebraic geometry \cite{FLPP}. 

The Weierstrass representation admits to apply to surface theory other branches of mathematics. 
In Section 2 we discuss an approach based on the spectral theory of the Dirac operator. 
The quaternionic algebraic geometry applies algebro-geometrical methods by considering solutions of the Dirac equation as ``holomorphic'' sections of spinor bundles.
It starts with treating the symmetry (\ref{symmetry}) as a multiplication by an imaginary unit
${\bf j}$  and considering $\mathrm{Ker} D$ as a linear space over quaternions ${\mathbb H}$ \cite{PP}. Therewith one may consider the Dirac operator of the more general form
\begin{equation}
\label{dirac4}
D = \left(\begin{array}{cc} 0 & \partial \\ -\bar{\partial} & 0 \end{array}\right) + 
\left(\begin{array}{cc} U & 0 \\ 0 & \bar{U} \end{array}\right)
\end{equation}
whose kernel is also invariant under (\ref{symmetry}).

For that we identify ${\mathbb C}^2$ with ${\mathbb H}$ as follows
$$
(z_1,z_2) \to z_1 + {\bf j} z_2 = \left(
\begin{array}{cc}
z_1 & -\bar{z}_2 \\
z_2 & \bar{z}_1
\end{array}
\right)
$$
and consider the  two matrix operators
$$
\bar{\partial} = \left(
\begin{array}{cc}
\bar{\partial} & 0 \\
0 & \partial
\end{array}
\right), \ \ \ {\bf j} U = {\bf j} \left(\begin{array}{cc} U & 0 \\ 0 &
\bar{U} \end{array}\right) = \left(\begin{array}{cc} 0 & -\bar{U} \\
U & 0 \end{array}\right),
$$
where ${\bf j} \in {\mathbb H}$ is the imaginary unit for which we have
${\bf j}^2 = -1, z{\bf j} = {\bf j} \bar{z}, \bar{\partial} {\bf j} = {\bf j} \partial$.
Then the Dirac equation $D\psi =0$ takes the form
$$
(\bar{\partial} + {\bf j} U)(\psi_1 + {\bf j} \psi_2) =
(\bar{\partial}\psi_1 - \bar{U}\psi_2) + {\bf j} (\partial \psi_2 + U \psi_1) = 0.
$$
Since $\psi_1$ and $\bar{\psi}_2$ are sections
of the same bundle $E$, we rewrite the Dirac equation as
$$
(\bar{\partial} + {\bf j} U)(\psi_1 + \bar{\psi}_2 {\bf j}) = 0
$$
and treat $E \oplus E$ as a quaternionic line bundle
whose sections are of the form $\psi_1 + \bar{\psi}_2 {\bf j}$. The symmetry (\ref{symmetry})
induces some quaternion linear endomorphism $J$ of $E$ such that
$J^2=-1$:  $\psi_1 + \bar{\psi}_2 {\bf j} \to (\psi_1 + \bar{\psi}_2 {\bf j}){\bf j} = -\bar{\psi_2} + \psi_1 {\bf j}$.
$J$ defines for any quaternion fiber a canonical splitting into ${\mathbb C} \oplus {\mathbb C}$ (in our case this is a splitting
into $\psi_1$ and $\bar{\psi}_2$) and such a bundle is called a ``complex quaternionic line bundle.''
The kernel of $D = \bar{\partial} + {\bf j}U$ is invariant under the right-side multiplications by constant quaternions and hence is a linear space over ${\mathbb H}$.

The ``quaternionic'' analog of the classical the Pl\"ucker formula established in \cite{FLPP}
 implies (\ref{conj1}) and (\ref{friedrich}).

By using the analogy with complex algebraic geometry, other interesting results were obtained, in particular on 
Backlund transformations and special classes of surfaces. Moreover this approach has another opportunity: in its framework the Weierstrass representation was also extended to surfaces in ${\mathbb R}^4$ and therewith ${\mathbb R}^4$ was naturally identified with ${\mathbb H}$. In the coordinate language the representation was written down in \cite{K2} and is as follows.

Let $D$ be of the form (\ref{dirac4}) and introduce the formally conjugate operator
$$
D^\vee = \left(\begin{array}{cc} 0 & \partial \\
-\bar{\partial} & 0 \end{array}\right) +
\left(\begin{array}{cc} \bar{U} & 0 \\
0 & U \end{array}\right).
$$

\begin{theorem}[\cite{K2}]
If $\psi$ and $\varphi$ satisfy the equations
\begin{equation}
\label{diract}
D\psi=0, \ \ \ \ D^\vee\varphi=0.
\end{equation}
then the formulae
$$
x^k(P) = x^k(P_0) + \int \left( x^k_z dz + \bar{x}^k_z
d\bar{z}\right), \ \ k=1,2,3,4,
$$ 
\begin{equation}
\label{weier4}
\begin{split}
x^1_z = \frac{i}{2} (\bar{\varphi}_2\bar{\psi}_2 + \varphi_1
\psi_1), \ \ \ \ x^2_z = \frac{1}{2} (\bar{\varphi}_2\bar{\psi}_2 - \varphi_1 \psi_1),
\\
x^3_z = \frac{1}{2} (\bar{\varphi}_2 \psi_1 + \varphi_1
\bar{\psi}_2), \ \ \ \ x^4_z = \frac{i}{2} (\bar{\varphi}_2 \psi_1 -
\varphi_1 \bar{\psi}_2),
\end{split}
\end{equation}
define the surface in ${\mathbb R}^4$ for which the induced metric is given by
$e^{2\alpha} dz d\bar{z} = (|\psi_1|^2+|\psi_2|^2)(|\varphi_1|^2+|\varphi_2|^2)dz d\bar{z}$
and $|U| = \frac{|{\bf H}| e^\alpha}{2}$ with ${\bf H}$ the mean curvature vector.
\end{theorem}

For $U=\bar{U}$ and $\psi = \varphi$ this representation reduces to (\ref{weier}).

The converse is also true but there is a difference with surfaces in ${\mathbb R}^3$ for which a choice of
a parameter $z$ defines $\psi$ uniquely up to multiplication by $\pm 1$.

\begin{theorem}[\cite{TDS}]
Every oriented surface (with a given conformal parameter) has representation (\ref{weier4}).
The spinors $\psi$ and $\varphi$ are defined up to the gauge transformations 
$$
\psi_1 \to e^h\psi_1, \ \psi_2 \to e^{\bar{h}}\psi_2, \ \varphi_1 \to e^{-h}\varphi_1, \ \varphi_2 \to e^{-\bar{h}}\varphi_2, \ U \to e^{\bar{h}-h}U,
$$
where $h$ is holomorphic. For every torus the potential $U$ may be taken double-periodic.
\end{theorem} 

Let us explain the appearance of these gauge transformations and, at the same time, why the dimensions 
$3$ and $4$ are distinguished by the existence of such spinor representations.

The Grassmannian $\widetilde{G}_{n,2}$ of oriented two-planes in ${\mathbb R}^n$ is diffeomorphic
to the quadric $Q$:
$$
z_1^2 + \dots + z_n^2 = 0, \ \ (z_1:\dots:z_n) \in Q_n  \subset {\mathbb C}P^{n-1}.
$$ 
To every oriented plane with an positively oriented orthonormal basis  
$e_1=(x_1,\dots,x_n), e_2 = (y_1,\dots,y_n)$ there corresponds the point 
$(z_1 : \dots : z_n), z_k = x_k + iy_k, k=1,\dots,n$, of this quadric.
Given a surface $(X^1(z,\bar{z}),\dots,X^n(z,\bar{z}))$ in ${\mathbb R}^n$ with a conformal parameter
$z$, we define the Gauss map as
$$
z \to \left(\frac{\partial X^1}{\partial z} : \dots : \frac{\partial X^n}{\partial z}\right) \in Q_n.
$$
It is straightforward  to derive that the image of the Gauss map lies in the quadric 
from the conformality of $z$.
For $n=3$ the quadric $Q_3$ is diffeomorphic to ${\mathbb C}^1$, its rational parameterization is
$$
z_1 = \frac{i}{2}(a^2 - b^2), \ z_2 = \frac{1}{2}(b^2-a^2), \ z_3 = ab, \ \ (a:b) \in {\mathbb C}P^1,
$$
and the spinor $\psi$ is reconstructed from the Gauss map: $\psi_1 = a, \bar{\psi}_2 = b$. 
For $n=4$ we have the diffeomorphic Segre mapping
$$
{\mathbb C}P^1 \times {\mathbb C}P^1 \to Q_4
$$ 
of the form $z_1 = \frac{i}{2}(a_1 b_1 + a_2 b_2), z_2 =
\frac{1}{2}(a_2 b_2 - a_1 b_1), z_3 = \frac{1}{2}(a_1 b_2 - a_2 b_1),  z_4 =
\frac{i}{2}(a_2 b_1 - a_1 b_2), (a_1: a_2) \in {\mathbb C}P^1, (b_1:b_2) \in {\mathbb C}P^1$.
The spinors take the form $\varphi = (a_1,\bar{a}_2), \psi = (b_1,\bar{b}_2)$ and are reconstructed up to 
the gauge transformations. Since they have to satisfy (\ref{diract}), $h$ has to be holomorphic.  
For $n >4$ the quadrics $Q_n$ have no such rational parameterizations.

The operators $D$  and $D^\vee$ enter  the representation of the Davey--Stewartson (DS) equations via compatibility of linear systems. That led to introducing the DS deformations of surfaces, the four-dimensional analog of the mNV deformations \cite{K2}. 

We consider one of such deformations for which we proved that it transforms tori into tori 
and preserves 
the Willmore functional $4\int |U|^2 dx \wedge dy$ \cite{TDS}. 
It has the form
\begin{equation}
\label{ds2}
U_t = i(U_{zz}+U_{\bar{z}\bar{z}} + (V+\bar{V})U), \ \ \ V_{\bar{z}} = 2(|U|^2)_z
\end{equation}
and is the compatibility condition for the linear problems
$$
D \psi =0 , \ \ \ \partial_t \psi = A\psi
$$
where
$$
A = i\left(\begin{array}{cc} -\partial^2 - V & \bar{U}\bar{\partial} - \bar{U}_{\bar{z}} \\
U\partial - U_z & \bar{\partial}^2 + \bar{V} \end{array}\right).
$$
It is also the compatibility condition for the system
$$
D^\vee \varphi = 0, \ \ \varphi_t = A^\vee \varphi,
$$
where
$$
A^\vee = - i\left(\begin{array}{cc} -\partial^2 - V & U\bar{\partial} - U_{\bar{z}} \\
\bar{U}\partial - \bar{U}_z & \bar{\partial}^2 + \bar{V} \end{array}\right).
$$
This equation is called the Davey--Ste\-wartson II (DSII) equation.

The evolution of $\psi$ and $\varphi$ gives us a deformation of the Gauss map of surfaces
(\ref{weier4}) which are at every moment of time defined up to a translation depending on the temporal variable. 

\section{The Moutard transformation for the Davey--Ste\-wartson II equation and its applications}

The Moutard transformation was  introduced in 1876 in projective differential geometry  for the equation 
$$
f_{xy} + Uf = 0.
$$
Given a solution $f_0$ of this equation, the transformation constructs another equation of this form with a different potential $\widetilde{U}$ such that to every solution of the first equation there corresponds a solution of the new one and
this is done by an explicit analytical formula. One of the problems to which the transformation was applied is an explicit construction of an immersion of the hyperbolic plane into ${\mathbb R}^3$ which, by Hilbert's theorem,  appeared to be impossible. Later the one-dimensional version, the Darboux transformation, was constructed and  has found many important applications in mathematical physics.  

Recently the version for the elliptic equation $f_{z\bar{z}} + Uf = 0$ was applied, for instance, to constructing in terms of explicit analytical formulae

1) blowing up solutions of the Novikov--Veselov equation with regular and fast decaying initial data \cite{TT},

2) two-dimensional von Neumann--Wigner potentials with multiple positive eigenvalues \cite{NTT}.

We recall that a potential of the Schr\"odinger operator on ${\mathbb R}^n$ is called von Neumann--Wigner if it has a positive eigenvalue.

Here we construct a Moutard type transformation for (\ref{diract}) and extend it to a transformation of solutions of the DSII equation.

Let extend spinors $\psi$ and $\varphi$ to ${\mathbb H}$-valued functions; i.e.,
$$
\Psi = \left(\begin{array}{cc} \psi_1 & -\bar{\psi}_2 \\
\psi_2 & \bar{\psi}_1 \end{array} \right), \ \ \ \
\Phi = \left(\begin{array}{cc} \varphi_1 & -\bar{\varphi}_2 \\
\varphi_2 & \bar{\varphi}_1 \end{array} \right)
$$
and put
$$
\omega(\Phi,\Psi) =
-\frac{i}{2}\left(\Phi^\top \sigma_3 \Psi + \Phi^\top \Psi\right) dz -
\frac{i}{2}\left(\Phi^\top \sigma_3 \Psi - \Phi^\top \Psi\right) d\bar{z},
$$
where $X \to X^\top$ is the conjugation of $X$, and
$\sigma_3= \left(\begin{array}{cc} 1 & 0 \\ 0 & -1 \end{array}\right)$ is the Pauli matrix.
If $\Psi$ and $\Phi$ satisfy the Dirac equations (\ref{diract}) then $\omega(\Phi,\Psi)$ and $\omega(\Psi,\Phi)$ 
are closed forms. Denote, for brevity, 
$\Gamma = \left(\begin{array}{cc} 0 & 1 \\ -1 & 0 \end{array}\right)$.
The ${\mathbb H}$-valued function
$$
S(\Phi,\Psi)(z,\bar{z}) = \Gamma \int \omega(\Phi,\Psi) = 
$$
$$
=
\int \left[ i\left(\begin{array}{cc} \psi_1\bar{\varphi}_2 & -\bar{\psi}_2\bar{\varphi}_2 \\
\psi_1 \varphi_1 & -\bar{\psi}_2\varphi_1 \end{array} \right) dz +
i\left(\begin{array}{cc} \psi_2\bar{\varphi}_1 & \bar{\psi}_1\bar{\varphi}_1 \\
- \psi_2 \varphi_2 & -\bar{\psi}_1\varphi_2 \end{array} \right) d\bar{z}\right]
 =
$$
$$
= \int d \left(\begin{array}{cc} ix^3 + x^4 & -x^1-ix^2 \\
x^1-ix^2 & -ix^3+x^4 \end{array}\right)
$$
defines a surface in ${\mathbb R}^4 = {\mathbb H}$ with $z$ the conformal parameter (\ref{weier4}).
Hence we identify $S$ with a surface in ${\mathbb R}^4$.

Let us define the ${\mathbb H}$-valued function
\begin{equation}
\label{kmatrix}
K(\Phi,\Psi) =  \Psi S^{-1}(\Phi,\Psi)\Gamma \Phi^\top\Gamma^{-1} =
\left(\begin{array}{cc} i\bar{W} & a \\ -\bar{a} & -iW \end{array}\right).
\end{equation}

The following theorem gives a Moutard type transformation for $D$.

\begin{theorem}[\cite{MT}]
Given $\Psi_0$  and $\Phi_0$, the solutions of (\ref{diract}),
for every pair $\Psi$ and $\Phi$  of solutions of the same equations
the ${\mathbb H}$-valued functions
$$
\widetilde{\Psi} =  \Psi - \Psi_0 S^{-1}(\Phi_0,\Psi_0) S(\Phi_0,\Psi), \
\widetilde{\Phi} =  \Phi - \Phi_0 S^{-1}(\Psi_0,\Phi_0) S(\Psi_0,\Phi)
$$
satisfy the Dirac equations
$$
\widetilde{D}\widetilde{\Psi} = 0, \ \ \ \ \widetilde{D}^\vee \widetilde{\Phi} = 0
$$
for the Dirac operators with the potential
\begin{equation}
\label{newpotential}
\widetilde{U} = U + W,
\end{equation}
where $W$ is defined by (\ref{kmatrix}) for $K(\Phi_0,\Psi_0)$.
Here $S(\Psi_0,\Phi_0)$ is normalized by the condition
$$
\Gamma S^{-1}(\Phi_0,\Psi_0)\Gamma = (S^{-1}(\Psi_0,\Phi_0))^\top
$$

The potential $\widetilde{U}$ is the potential of the Weierstrass representation of the surface $S^{-1}$
with $z$ a conformal parameter. The surface $S^{-1}$ is obtained from $S$ by composition of 
the inversion centered at the origin and the reflection
$(x_1,x_2,x_3,x_4) \to (-x_1,-x_2,-x_3,x_4)$.
\end{theorem}

For $U = \bar{U}$ and $\Psi = \Phi$  this transformation reduces to the transformation of Dirac operators with real-valued potentials given in \cite{C} in different form. In \cite{T151} it was related to the Weierstrass representation of surfaces in ${\mathbb R}^3$ by proving that it corresponds to the Mobius inversion $S \to S^{-1}$ 

Let us replace $K(\Phi,\Psi)$ in (\ref{kmatrix}) with
$$
S(\Phi,\Psi)(z,\bar{z},t) = \Gamma \int \omega(\Phi,\Psi) + \Gamma \int \omega_1(\Phi,\Psi),
$$
where
$$
\omega_1(\Phi,\Psi) = 
\left(\left[
\Phi^\top_z \left(\begin{array}{cc} 1 & 0 \\ 0 & 0 \end{array}\right) +
\Phi^\top_{\bar{z}}\left(\begin{array}{cc} 0 & 0 \\ 0 & 1 \end{array}\right)
\right] \Psi
\right.
$$
$$  
\left.
- \Phi^\top\left[
\left(\begin{array}{cc} 1 & 0 \\ 0 & 0 \end{array}\right)\Psi_z +
\left(\begin{array}{cc} 0 & 0 \\ 0 & 1 \end{array}\right)\Psi_{\bar{z}}\right]\right) dt.
$$
We have

\begin{theorem}[\cite{T2021}]
If $U$ meets the Davey--Stewartson II equation (\ref{ds2}) and $\Psi$ and $\Phi$ satisfy the equations
$D \Psi = 0, \Psi_t = A\Psi, D^\vee \Phi =0, \Phi_t = A^\vee \Phi$, then the Moutard transformation (\ref{newpotential}) of $U$ gives the solution $\widetilde{U}$ of the DSII equation
$$
\widetilde{U}_t = i(\widetilde{U}_{zz}+\widetilde{U}_{\bar{z}\bar{z}} + 2(\widetilde{V}+\bar{\tilde{V}})\widetilde{U}), \ \ \
\widetilde{V}_{\bar{z}} = (|\widetilde{U}|^2)_z
$$
with
\begin{equation}
\label{newv}
\widetilde{V} = V + 2ia_z
\end{equation}
where $a$ is given by (\ref{kmatrix}).
\end{theorem}

The geometrical meaning of this transformation is as follows:
for every fixed $t$ the spinors $\Psi$ and $\Phi$ determine some surface $S(t)$ in ${\mathbb R}^4$  and $U$ is the potential of such a representation. The surfaces $S(t)$ evolve via the DSII equation.
We invert every such a surface and obtain the $t$-parameter family of surfaces $\widetilde{S}(t) = S^{-1}(t)$ which evolve via the DSII equation.
Starting with a family of smooth surfaces and the corresponding smooth potentials $U$ we may construct singular solutions of the DSII equation: when $S(t)$ passes through the origin the function $\widetilde{U}$ loses continuity or regularity because the origin is mapped into the infinity by the inversion.

One of the simplest applications of Theorem 1 consists in constructing exact solutions 
from holomorphic functions. In this case we start from the trivial solution $U=V=0$ for which $\Psi$ and $\Phi$ are defined by holomorphic data. For instance, we have

\begin{theorem}[\cite{T2021}]
\label{lasttheo}
Let $f(z,t)$ be a function which is holomorphic in $z$ and satisfies the equation
$$
\frac{\partial f}{\partial t} = i\frac{\partial^2 f}{\partial z^2}.
$$
Then
$$
U = \frac{i(zf^\prime - f)}{|z|^2 + |f|^2}, \ \
V = 2ia_z,
$$
where
$$
a = -\frac{i(\bar{z} + f^\prime) \bar{f}}{|z|^2 + |f|^2}
$$
satisfy the Davey--Stewartson II equation.
\end{theorem}

Geometrically we have the deformation of graphs $w=f(z,t)$ which are minimal surfaces in ${\mathbb R}^4 = {\mathbb C}^2$. Whenever $f(z,t)$ vanishes at $z=0$ the graph passes through the origin and the solution $\widetilde{U}$ loses continuity or regularity. Hence the Weierstrass representation visualizes the creation of singularity and gives a method for finding such solutions. 

We already applied this idea to constructing a solution with a one-point singularity for the modified Novikov--Veselov equation by using the Enneper surface \cite{T152}. However in difference with the mNV equation the DSII has an important physical meaning.

In the variables
$$
X = 2y, \ \ Y = 2x.
$$
the Davey--Stewartson II equation takes the form known in mathematical physics, namely,
\begin{equation}
\label{ozeq}
iU_t -U_{XX} + U_{YY} = -4|U|^2U + 8\varphi_X U,
\end{equation}
$$
\Delta\varphi = \frac{\partial^2 \varphi}{\partial X^2}
+ \frac{\partial^2 \varphi}{\partial Y^2} = \frac{\partial}{\partial X} |U|^2,
$$
where $\mathrm{Re} V = 2|U|^2 - 4 \varphi_X$, $\varphi_X = \frac{\partial \varphi}{\partial X}$ \cite{DS}.
This version of the DSII equation is called focusing.

Ozawa constructed a blow-up solution to (\ref{ozeq})  with the initial data
$$
U(X,Y,0) = \frac{e^{-ib(4a)^{-1}(X^2 - Y^2)}}{a(1+ ((X/a)^2+(Y/a)^2)/2)}
$$
and showed that for constants $a$ and $b$ such that $ab<0$ we have
$$
\| U\|^2 \to 2 \pi\cdot \delta \ \mbox{as $t \to T = -a/b$}
$$
in ${\cal S}^\prime$ where 
$$
\|U\|^2 = \int_{{\mathbb R}^2}|U|^2\, dx\,dy
$$ 
is the squared $L_2$-norm of $U$ and $\delta$ is the Dirac distribution centered at the origin. We remark that $\|U\|^2 = 2\pi$ and the solution extends to $T > -a/b$ and gains regularity.
In \cite{KS} it is conjectured for this equation that the blow-up in all cases is self-similar and  the time dependent 
scaling is as in the Ozawa solution. This conjecture is based on numerical results.

Let us consider the simplest examples of the solutions given by Theorem \ref{lasttheo}.
We denote  by $c$ a constant $c$ which may take arbitrary complex values and by $r$ we
denote $|z|, z \in {\mathbb C}$.

1) $f = z^2 + 2it + c$: 
\begin{equation}
\label{s1}
U = \frac{i(z^2-2it-c)}{|z|^2 + |z^2 + 2it+c|^2}, \
\end{equation}
$$
V =
\frac{4(\bar{z}^2-2it + \bar{c})}{|z|^2 + |z^2+2it+c|^2}-
\frac{2(2z(\bar{z}^2 -2it +\bar{c})+\bar{z})^2}{(|z|^2 + |z^2+2it+c|^2)^2},
$$
$|U| = O\left(\frac{1}{r^2}\right)$ as $r \to \infty$.
If $c$ is not purely imaginary, then the solution is always smooth.
It $c = it\tau, \tau \in {\mathbb R}$, then for $t = -\frac{\tau}{2}$  $U$ has singularity at $z=0$ of the type
$$
U \sim i e^{2i \phi} \ \ \ \mbox{as $r  \to 0$}, \ \ \mbox{where $z = r e^{i\phi}$}.
$$
Remark that $U \in L_2({\mathbb R}^2)$ for all $t$ and $c$.
Since a small variation of $c$ removes singularities, they are unstable.

2) $f=z^4+12itz^2 - 12t^2 +c$,
\begin{equation}
\label{s2}
U = \frac{i(3z^4 + 12 itz^2 + 12t^2 -c)}{|z|^2 + |z^4 +12itz^2 - 12t^2 +c|^2}.
\end{equation}
This solution becomes singular for $c = 12t^2$ which is possible if and only if
$c$ is real-valued and positive. In this case it has singularities $U \sim -12t e^{2i\phi}$ at $z=0$ for $t=\pm \sqrt{c/12}$.

We remark that $\|U\|^2$ is the first integral of the system. 
For (\ref{s1}) it is always equal to $2\pi$ except for the time $T_{\mathrm{sing}}$ when the solution becomes singular. For $t= T_{\mathrm{sing}}$ it is equal to
$\pi$. Analogously for (\ref{s2}) it  is equal to $4\pi$ for $t$ such that $U$ is nonsingular and is equal to $3\pi$ for $t=T_{\mathrm{sing}}$. The multiplicity of the value of this functional to $\pi$
in both cases is explained by that the surfaces $\widetilde{S}$ are immersed Willmore spheres (with singularities for singular moments of time).

By taking polynomials of higher degrees for $f$, we can construct such singular solutions 
for which the regular initial data have any polynomial decay. 

Are there another physically relevant wave equations that admit solutions with such singularities?

\end{document}